\documentclass{amsart}
\usepackage{amsfonts, amsmath, latexsym, epsfig}
\usepackage[norelsize]{algorithm2e}
\usepackage{amssymb}
\usepackage{color}
\usepackage{epsf}
\usepackage{url}

\newtheorem{theorem}{Theorem}

\newtheorem{conjecture}{Conjecture}

\def\QuotS#1#2{\leavevmode\kern-.0em\raise.2ex\hbox{$#1$}\kern-.1em/\kern-.1em\lower.25ex\hbox{$#2$}}

\newcommand{\MF}{{\mathcal F}}

\DeclareMathOperator{\Aut}{Aut}

\DeclareMathOperator{\CUT}{CUT}
\DeclareMathOperator{\CUTP}{CUTP}

\DeclareMathOperator{\CORP}{CORP}
\DeclareMathOperator{\MET}{MET}
\DeclareMathOperator{\METP}{METP}

\begin{document}

\author{Michel Deza}
\address{Michel Deza, \'Ecole Normale Sup\'erieure, 75005 Paris, France}
\email{Michel.Deza@ens.fr}

\author{Mathieu Dutour Sikiri\'c}
\address{Mathieu Dutour Sikiri\'c, Rudjer Boskovi\'c Institute, Bijenicka 54, 10000 Zagreb, Croatia, Fax: +385-1-468-0245}
\email{mathieu.dutour@gmail.com}

\title{Enumeration of the facets of cut polytopes over some highly symmetric graphs}

\date{}

\maketitle

\begin{abstract}
We report here a computation giving the complete list of facets for the cut polytopes
over several very symmetric graphs with $15-30$ edges, including $K_8$, $K_{3,3,3}$,
$K_{1,4,4}$, $K_{5,5}$, some other $K_{l,m}$, $K_{1,l,m}$,
$Prism_7, APrism_6$, M\"{o}bius ladder $M_{14}$, Dodecahedron, Heawood and Petersen graphs.

For $K_8$, it shows that the huge lists of facets of the cut polytope $\CUTP_8$ and
cut cone $\CUT_8$, given in \cite{CR} is complete. We also confirm the conjecture that
any facet of $\CUTP_8$ is adjacent to a triangle facet.

The lists of facets for $K_{1,l,m}$ with $(l,m)=(4,4),(3,5),(3,4)$ solve problems (see,
for example, \cite{WernerBellInequality}) in quantum information theory.
\end{abstract}

\section{Introduction}
Polyhedra, generated by cuts and by finite metrics are central objects of Discrete Mathematics; see, say, \cite{DL}.
In particular, they are tightly connected with the well-known  NP-hard  optimization problems
such as the max-cut problem and the unconstrained quadratic $0,1$ programming problem.
To find their (mostly unknown) facets is one approach to this problem.
It is a difficult problem since the number of facets grows very fast with the size of the problem.
However, when full computation is unfeasible, having a partial list of facets can help branch-and-bound
strategies that are then commonly used. Here we consider the case of cut polytope over graphs.

Given a graph $G=(V,E)$ with $n=|V|$, for a vertex subset $S\subseteq V=\{1, \dots, n\}$, the {\em cut semimetric} $\delta_S(G)$ is a vector (actually, a symmetric $\{0,1\}$-matrix) defined as
\begin{equation*}
\delta_S(x,y)=\left\{\begin{array}{rcl}
1 &\mbox{~~~if~~~} x y\in E \mbox{~~~~and ~~~}  \vert S \cap \{x,y\}\vert = 1\\
0 &\mbox{otherwise.}
\end{array}\right.
\end{equation*}
So, $\delta_S$ can be seen also as the adjacency matrix of a {\em cut} (into $S$ and $\overline{S}$) {\em subgraph} of $G$. 
Clearly, $\delta_{\{1, \dots, n\} - S} = \delta_S$.
If $G$ is connected, which will be the case in our work, there are exactly $2^{n-1}$ distinct cut semimetrics.
The {\em cut polytope} $\CUTP(G)$ and the {\em cut cone} $\CUT(G)$ are defined as the convex hull of all such  semimetrics and the positive span of all non-zero ones among them, respectively. Their number of vertices, respectively extreme rays is $2^{n-1}$, respectively $2^{n-1} - 1$ and their dimension is $\vert E\vert$, i.e., the number of edges of $G$.

The most interesting and complicated case is $\CUTP(K_n)$ and $\CUT(K_n)$, denoted simply
$\CUTP_n$, $\CUT_n$ and called {\em the cut polytope} and  {\em the cut cone}.
In fact, $\CUT_n$ is the set of all $n$-vertex semimetrics, which embed isometrically into some
metric space $l_1$, and rational-valued elements of $\CUT_n$ correspond exactly to the $n$-vertex
semimetrics, which embed isometrically, {\em up to a scale} $\lambda \in \mathbb{N}$, into 
the path metric of some $N$-cube $K_2^N$. It shows importance of this cone in Analysis and
Combinatorics.

The {\em metric cone} $\MET_n$ is  the set of all {\em semimetrics on $n$ points}, i.e., those of above functions,
which  satisfy all {\em triangle inequalities}, i.e., $d_{ij}\leq d_{ik} + d_{kj}$.
The bounding of $\MET_n$ by ${n\choose 3}$ {\em perimeter inequalities} 
$d_{i j}+d_{i k}+d_{j k} \leq 2$ produces the {\em metric polytope} $\METP_n$. 
We have the evident inclusions $\CUT_n \subseteq \MET_n$ and $\CUTP_n  \subseteq \METP_n$ with 
$\CUT_n = \MET_n$ and $\CUTP_n = \METP_n$ only for $3 \le n \le 4$. Another relaxation of 
the cut-polytope is the {\em hypermetric polytope} considered in \cite{Hyp8}.
In Table \ref{ResultsCUTn_METn} we give the number of facets, vertices (or extreme rays) and
the number of orbits for $\CUTP_n$, $\METP_n$ with  $n\leq 8$ and for the corresponding cones.

The symmetry group of a graph $G=(V,E)$ induces symmetry of $\CUTP(G)$.
For any $U\subset \{1,\dots, n\}$, the map $\delta_S\mapsto \delta_{U\Delta S}$ also defines
a symmetry of $\CUTP(G)$. Together those form the {\em restricted symmetry group} $ARes(\CUTP(G))$,
the order of which is $2^{|V|-1}\vert \Aut(G)\vert$. The full symmetry group $\Aut(\CUTP(G))$
may be larger.
For $n\not= 4$, $ \Aut(\CUTP(K_n))=ARes(\CUTP(K_n))$ but $\Aut(\CUTP(K_4))=$ $\Aut(K_{4,4})$ (\cite{Relatives}).
So, $|\Aut(\CUTP(K_n))|$ is $2^{n-1}n!$ for $n\neq 4$ and $6\times 2^34!$ for $n=4$.
Let $A(G)$ denote $2^{1-|V|}|\Aut(\CUT(G))|$. 
One can check the following.
\begin{enumerate}

\item If $G$ is a complete multipartite graph with $t_i\ge 1$ parts of size $a_i$ for $1\le i \le r$ and $1\le a_1< a_2< \dots <a_r$, then $|Aut(G)|=\prod_{i=1}^rt_i!(a_i!)^{t_i}$.

\item If $G$ is  $Prism_m$, $APrism_m$ or  M\"{o}bius ladder $M_{2m}$, then
$A(G)=|Aut(G)|$ and,  moreover, $A(G)=4m=2|V|$  if $m\neq 4$, $m\ge 4$  and 
$m\ge 4$, respectively.

\item Examples of infinite series of {\em singular}, i.e., having
$|Aut(\CUTP(G))|>|ARes(\CUTP(G))|$, graphs $G$ are $K_{2,m}$, $K_{1,1,m}$
for any $m\ge 2$ and $P_m$, $C_m$  for any $m\ge 3$.
All $4$-, $5$- and $6$-vertex connected graphs are singular, except 
$K_{1,3}$, $K_{1,4}$,  $K_5$, $K_5-e$, $K_5-e_1-e_2$ with disjoint
edges $e_1,e_2$ and $19$ (among all $112$) $6$-vertex graphs.
\end{enumerate}

In Table \ref{TableSomeCUTpolytopes} we list information on the cut polytopes of several graphs.
Full data sets are available from \cite{WebPageCutPolytopes}. The computational techniques used
are explained in Section \ref{SEC_Computational_Methods}; they can be applied to any polytope having
a large symmetry group.
In Section \ref{SEC_CUTP8} the results for $\CUTP_8$ are detailed.
In Section \ref{SEC_CorrelationPolytope} the results for the correlation polytopes $K_{n,m}$
are presented.
Finally, in Section \ref{SEC_GraphsWithoutK5minor} we explain the results for graphs without $K_5$ minors.

\begin{table}
\begin{center}
\caption{The number of facets and vertices (or extreme rays) in  the cut and metric polytopes (or cones) for $n\le 8$.
The enumeration of orbits of facets of $\CUT_n$ and $\CUTP_n$ for $n\le 7$ was done in
\cite{OrbitFacetCutPolytope5,OrbitFacetCutPolytope6,OrbitFacetCutPolytope7} for $n=5$, $6$
and $7$, respectively.
The enumeration of orbits of extreme rays of $\MET_7$ was done in \cite{OrbitRaysMetricCone7}.
The orbits of vertices of $\METP_n$ were enumerated in \cite{OrbitVerticesMetricPolytope7}
for $n=7$ and in \cite{OrbitPolytopeMetricPolytope8} for $n=8$. For $n\le 6$ such enumeration
is easy.}
\label{ResultsCUTn_METn}
\begin{tabular}{||c||c|c|c|c|c|c||}
\hline
\hline
$P$               & $n=3$  & $n=4$   & $n=5$   & $n=6$    & $n=7$         & $n=8$\\ \hline\hline
$\CUTP_n,e$       & $4(1)$ & $8(1)$  & $16(1)$ & $32(1)$  & $64(1)$       & $128(1)$\\
$\bf {\CUTP_n,f}$ & $4(1)$ & $16(1)$ & $56(2)$ & $368(3)$ & $116,764(11)$ & $\bf{217,093,472(147)} $\\\hline
$\CUT_n, e$       & $3(1)$ & $7(2)$  & $15(2)$ & $31(3)$  & $63(3)$       & $127(4)$\\
$\bf {\CUT_n,f}$  & $3(1)$ & $12(1)$ & $40(2)$ & $210(4)$ & $38,780(36)$  & $\bf{49,604,520(2,169)}$\\
\hline
$\MET_n,e$        & $3(1)$ & $7(2)$  & $25(3)$ & $296(7)$ & $55,226(46)$  & $119,269,588(3,918)$\\
$\MET_n,f$        & $3(1)$ & $12(1)$ & $30(1)$ & $60(1)$  & $105(1)$      & $168(1)$\\
\hline
$\METP_n,e$       & $4(1)$ & $8(1)$  & $32(2)$ & $554(3)$ & $275,840(13)$ & $1,550,825,600(533)$\\
$\METP_n,f$       & $4(1)$ & $16(1)$ & $40(1)$ & $80(1)$  & $140(1)$      & $224(1)$\\
\hline
\hline\end{tabular}
\end{center}
\end{table}

\begin{table}
\begin{center}
\caption{The number of facets of the cut polytopes $\CUTP(G)$ of some graphs.
In $4$-th column, $A(G)$ denote $2^{1-|V|}|\Aut(\CUTP(G))|$. 
The case $A(G) >|\Aut(G)|$ is indicated by ${}^{*}$} 
\label{TableSomeCUTpolytopes}
\begin{tabular}{||c|c|c||c|c||}
\hline\hline
$G=(V,E)$             & $\vert V\vert$ & $\vert E\vert$ & $A(G)$ & Number of facets (orbits)\\\hline\hline
$K_8$                 & $8$   & $28$   & $8!$ & $\bf{217,093,472(147)}$\\\hline
$K_{3,3,3}$           & $9$   & $27$   & $(3)^4$ & $\bf{624,406,788(2015)}$\\
$K_{1,4,4}$           & $9$   & $24$   & $2(4!)^2$ & $\bf{36,391,264(175)}$\\
$K_{1,3,5}$           & $9$   & $23$   & $3!5!$ & $\bf{71,340(7)}$\\
$K_{1,3,4}$           & $8$   & $19$   & $3!4!$ & $\bf{12,480(6)}$\\
$K_{1,3,3}$           & $7$   & $15$   & $2(3!)^2$ & $\bf{684(3)}$\\
$K_{1,1,3,3}$         & $8$   & $21$   & $4(3!)^2$ & $\bf{432,552(50)}$\\
$K_{1,2,m}$, $m\ge 2$ & $m+3$ & $3m+2$ & $|\Aut(K_{1,2,m})|$  &  $\bf{8m + 8{m\choose 2}(2)}$\\
\hline
$K_{5,5}$             & $10$  & $25$   & $2(5)^2$ & $\bf{16,482,678,610(1,282)}$\\
$K_{4,7}$             & $11$  & $28$   & $4!7!$ & $\bf{271,596,584(15)}$\\
$K_{4,6}$             & $10$  & $24$   & $4!6!$ & $\bf{23,179,008(12)}$\\
$K_{4,5}$             & $9$   & $20$   & $4!5!$ & $\bf{983,560(8)}$\\
$K_{4,4}$             & $8$   & $16$   & $2(4!)^2$ & $\bf{27,968(4)}$\\
$K_{3,m}, m\ge 3$     & $m+3$ & $3m$   & $|\Aut(K_{3,m})|$ & $\bf{6m + 24{m\choose 2}(2)}$\\
$K_{2,m},m\ge 3$      & $m+2$ & $2m$   & $2^{m-1}m!|\Aut(K_{2,m})|$ $\,^{*}$  & $\bf{4m^2(1)}$\\
$K_{1,m},m\ge 2$      & $m+1$ & $m$    & $m!$   & $\bf{2m(1)}$\\
$K_{m+2}-K_m, m\ge 2$ & $m+2$ & $2m+1$ & $2^{m-1}m!|\Aut(K_{1,1,m})|$ $\,^{*}$  & $\bf{4m(1)}$\\
$K_{m+3}-K_m, m\ge 2$ & $m+3$ & $3m+3$ & $3!m!$ & $\bf{4+12m(2)}$\\
$K_{m+4}-K_m, m\ge 2$ & $m+4$ & $4m+6$ & $4!m!$ & $\bf{8(8m^2-3m+2)(4)}$\\
$K_8-K_3$             & $8$   & $25$ & $360$ & $\bf{2,685,152(82)}$\\
$K_7-K_2$             & $7$   & $20$ & $240$ & $\bf{31,400(17)}$\\
\hline
Dodecahedron          & $20$  & $30$ & $120$ & $\bf{23,804(5)}$\\
Icosahedron           & $12$  & $30$ & $120$ & $\bf{1,552(4)}$\\
Cube                  & $8$   & $12$ & $48$ & $\bf{200(3)}$\\
Cuboctahedron         & $12$  & $24$ & $48$ & $\bf{1,360(5)}$\\
Tr. Tetrahedron       & $12$  & $18$ & $24$ & $\bf{540(4)}$\\
$APrism_6$            & $12$  & $24$ & $24$ & $\bf{2,032(5)}$\\
$Prism_7$             & $14$  & $21$ & $28$ & $\bf{7,394(6)  }$\\
\hline
$Pyr(Prism_5)$        & $11$  & $25$ & $20$ & $\bf{208,132(22)  }$\\
$Pyr(APrism_4)$       & $9$   & $24$ & $16$ & $\bf{389,104(17)  }$\\
M\"{o}bius ladder $M_{14}$ & $14$ & $21$              & $28$    & $\bf 369,506{(9)}$\\
Heawood graph         & $14$  & $21$ & $336$ & $\bf{5,361,194(9)}$\\
Petersen graph        & $10$  & $15$ & $120$ & $\bf{3,614(4)}$\\
\hline\hline
\end{tabular}
\end{center}
\end{table}

\section{Computational methods}\label{SEC_Computational_Methods}
The second author has developed over the years an effective computer program (\cite{Polyhedral})
for enumerating facets of polytopes which are symmetric.
The technique used is {\em adjacency decomposition method},
 originally introduced in \cite{CR} and applied to the Transporting Salesman polytope,
the Linear Ordering polytope and the cut polytope.
 The algorithm is detailed in Algorithm \ref{AdjacencyDecompositionMethod} and surveyed in \cite{symsurvey}.

The initial facet of the polytope
$P$ is obtained via linear programming. The tests of equivalence are done via the {\tt GAP}
functionality of permutation group and their implementation of partition backtrack.
The problematic aspect is computing the facets adjacent to a facet. This is itself a dual description
problem for the polytope defined by the facet $F$.

\begin{figure}
\centering
\begin{minipage}[b]{5.5cm}
\centering
\epsfig{file=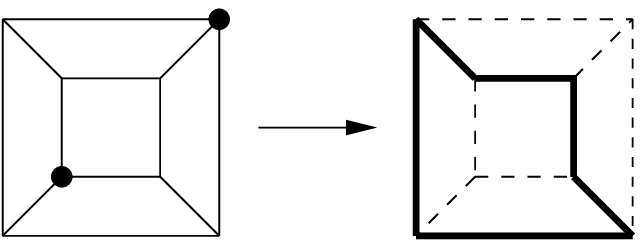,width=5.5cm}
a) A cube and the remaining connected graph obtained from Theorem \ref{ThConnectivityTheorem}.(i)
\end{minipage}
\begin{minipage}[b]{5.5cm}
\centering
\epsfig{figure=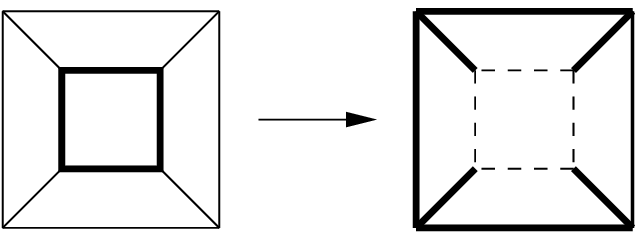,width=5.5cm}
b) A cube and the remaining connected graph obtained from Theorem \ref{ThConnectivityTheorem}.(ii)
\end{minipage}
\caption{Illustration of connectivity results of Theorem \ref{ThConnectivityTheorem}.}
\label{FigConnectivityTheorem}
\end{figure}

An interesting problem is the check when ${\mathcal R}$ is complete.
Of course, if all the orbits are treated, i.e., if ${\mathcal R}={\mathcal D}$, then the computation
is complete.
However, sometimes we can conclude before that:
\begin{theorem}\label{ThConnectivityTheorem}
Let $G(P)$ be the skeleton graph of a $m$-dimensional polytope.

(i) (\cite{BalinskiTh}) $G(P)$ is at least $m$-connected.

(ii) If we remove all the edges contained in a given face $F$, then the remaining graph is still connected.
\end{theorem}
Hence, if the total set of facets, equivalent to a facet in ${\mathcal R}\diagdown {\mathcal D}$, has
size at most $m-1$ or contains a common vertex, then we can conclude that ${\mathcal R}$ is complete.
The requirements for applying any of those two criterions are rather severe, but they are easy
to check and, when applicable, the benefits are large.
The geometry underlying Theorem \ref{ThConnectivityTheorem} is illustrated in Figure \ref{FigConnectivityTheorem}.

\begin{algorithm}\label{AdjacencyDecompositionMethod}
\KwData{Polytope $P$ and a group $G$}
\KwResult{Set~${\mathcal R}$ of all inequivalent representative of facets of $P$ for $G$}
$F \leftarrow$ facet of $P$.\\
${\mathcal R} \leftarrow \{ F \}$.\\
${\mathcal D} \leftarrow \emptyset$\\
\While{${\mathcal R}$ is not complete}{
  $F$ a facet in ${\mathcal R}\diagdown {\mathcal D}$.\\
  ${\mathcal D}\leftarrow {\mathcal D} \cup \{F\}$.\\
  $\MF \leftarrow$ facets of $P$ adjacent to $F$.\\
  \For{$F \in \MF$}{
    ${\rm test} \leftarrow {\rm true}$\\
    \If{$F$ is not equivalent to a facet in ${\mathcal R}$}{
      ${\mathcal R}\leftarrow {\mathcal R} \cup \{F\}$.\\
    }
  }
}
\caption{The adjacency decomposition method}
\end{algorithm}

The adjacency decomposition method works well when the polytope is symmetric but it still
relies on computing a dual description. This is easy when the incidence of the facet is low
but become more and more problematic for facets of large incidence. However, quite often 
high incidence facets also have large symmetry groups. Therefore, a natural extension of
the technique is to apply
the method recursively and so, obtain the {\em recursive adjacency decomposition method}, which
is again surveyed in \cite{symsurvey}.

In order to work correctly, the method requires several ingredients. One is the ability to
compute easily automorphism group of the polytope, see \cite{LMS_JCM_groupPolytope} for details.
Another is good heuristics for deciding when to apply the method recursively or not.
Yet another is a storing system for keeping dual description that may be reused, and this again depends
on some heuristics.
The method has been applied successfully on numerous problems
\cite{PerfectDim8,ComplexityVoronoiDSV,ContactLeech} and here.

The framework, that we have defined above, can be applied, in order
to sample facets of a polytope, say,  $P$. A workable idea is to use linear programming
and then get some facets of $P$. But doing so, we overwhelmingly get
facets of high incidence, while we may be interested in obtaining facets of low incidence.
One  can adapt the adjacency decomposition method to do such a sampling. Let us call
two facets {\em equivalent} if their incidence is the same. By doing so, we remove the combinatorial
explosion, which  is the main difficulty of such dual-description problem.
At the end, we get a number of orbits of facets of different incidence, which give an idea
of the complexity of the polytope.

\section{The facets of $\CUTP_8$}\label{SEC_CUTP8}
In \cite{CR}, the adjacency decomposition method was introduced and was applied to the 
Transporting Salesman polytope, the Linear Ordering polytope and the cut polytope.
For the cut polytope $\CUTP_8$, the authors found $147$ orbits, consisting of $217,093,472$ facets,
 but this
list was potentially incomplete, since their method did not treat
the triangle,
pentagonal and $7$-gonal inequalities (defined in Section 3) at that time.
Therefore, they only prove that the number of orbits is at least $147$.
The enumeration is used in \cite{AvisCorr,AvisTwoParty} for work in quantum mechanics.

Sometimes (\cite{DL,PolymakeLin}), this is incorrectly understood and it is reported that the
number of orbits is exactly $147$ with \cite{CR} as a reference.
Here  we show that Christof-Reinelt's list is complete. We had to treat the $3$ remaining
orbits of facets since we could not apply Theorem \ref{ThConnectivityTheorem}. However, we
could apply the theorem in deeper levels of the recursive adjacency method and this made
the enumeration faster.

\begin{conjecture}(\cite{DualCutPolytope,DL})
Any facet of $\CUTP_n$ is adjacent to a triangle inequality facet.
\end{conjecture}
The conjecture was checked for $n\leq 7$; here we confirm it for $n=8$.

Looking for a counterexample to this conjecture, we applied our sampling framework
to $\CUTP_n$ for $n=10$, $11$ and $12$. We got initial facets of low incidence and then we 
complemented this with random walks in the set of all facets. This allowed us to find many
simplicial facets (more than $10,000$ for each) of these $\CUTP_n$ but all of them were adjacent
to at least one triangle inequality facet.

As a corollary to the confirmation of above conjecture for $n=8$, we obtain that the ridge graphs of $\CUTP_8$ and $\CUT_8$   have diameter $4$.
In fact, the subgraphs, formed by the triangle facets, have diameter $2$.

\section{Correlation polytopes of $K_{n,m}$}\label{SEC_CorrelationPolytope}

In quantum physics and quantum information theory, 
 Bell inequalities, involving joint probabilities of two probabilistic events, are exactly inequalities valid for 
 the {\em correlation polytope} $\CORP(G)$   (called also {\em Boolean quadric polytope} $BQP(G)$)  
 of a graph, say, $G$.
 In particular, $\CORP(K_{n,m})$
   is seen in quantum theory as the set of possible results of a
 series of Bell experiments with a non-entangled (separable) quantum state shared by two distant parties,   where one party (Alice) has $n$ choices of possible two-valued measurements and the other party (Bob) has $m$ choices.  
 Here, a valid inequality of $\CORP(K_{n,m})$ is called a {\em Bell inequality} and if facet inducing, a {\em tight Bell inequality}. 
This polytope is equivalent (linearly isomorphic via the {\em covariance map}) to the cut polytope of $K_{1,n,m}$ \cite[Section 5.2]{DL}. Similarly,  $\CUTP(K_{1,n,m,l})$ represents three-party  Bell inequalities.

The symmetry group of  $\CORP(K_{n,n})$ and  $\CUTP(K_{1,n,m})$  is of order $2^{1+n+m}n!m!$.

We computed the 
 facets of $\CORP(K_{n,m})$ having $(n,m)=(2,m)$ with $m\le 6$, $(3,4),(3,5),(4,4)$ and confirmed known enumerations for $(n,m)=(2,2),(3,3)$. 
  In fact, the page \cite{WernerBellInequality}
 collects the progress 
 on finding Bell inequalities.
 The case of $\CORP(K_{n,n})$ is
called   there {\em $(2,n,2)$-setting}. 
The cases $n=2,3$ were settled in \cite{FineHidden} and \cite{PitowskyOptimal}, respectively.
For $n=4$, 
 partial lists of facets were known; 
 our $175$ orbits of $36,391,264$ facets of $\CORP(K_{4,4})$
 finalize this case.

In contrast to the Bell's inequalities, which probe entanglement
between spatially-separated systems, the {\em Leggett-Garg inequalities}
test the correlations of a single system measured at different times. The
polytope, defined by those inequalities  for $n$ observables,
is, actually (\cite{LeggettGargIneq}), the cut polytope $\CUTP_n$.

\section{Cut polytopes of graphs without minor $K_5$}\label{SEC_GraphsWithoutK5minor}

Given a graph $G=(V,E)$, an edge $(v_1v_2)\in E$ defines an {\em edge inequality}
$x(v_1, v_2)\geq 0$. 
Similarly, an $s$-cycle $(v_1,\dots , v_s)$ of $G$ with $s\geq 3$ defines a
{\em cycle inequality}
\begin{equation*}
\sum_{i=1}^{s-1} x(v_i,v_{i+1})  -x(v_1,v_s)\ge 0.
\end{equation*}
These inequalities and their switching define valid inequalities on $\CUTP(G)$.
There are $2|E|$ edge inequalities and their incidence is $2^{|V|-2}$.
Each $s$-cycle of $G$ gives $2^{s-1}$ $s$-cycle faces, which are of incidence $s2^{|V|-s}$.

In fact, \cite{BarahonaMahjoubOnCutPolytope,BarahonaCutContract} and \cite{SeymourMatroidMulticommodity}, gives, respectively, that:

(i)  only chordless $s$-cycles and edges not containing in a $3$-cycle produces facets;

(ii) no other facets exist if and only if $G$ is a $K_5$-minor-free graph.
 
All $K_5$-minor-free graphs in Table \ref{TableSomeCUTpolytopes} are $K_{1,2,m}$; $K_{2,m}$; $K_{m+i}-K_{m}$
with $m\ge 2, 1\le i\le 3$ and the skeletons of regular and semiregular polyhedra (Dodecahedron, Icosahedron,
Cube, Cuboctahedron, truncated Tetrahedron, $APrism_6, Prism_7$).

\section{Acknowledgments}

Second author gratefully acknowledges support from the Alexander von Humboldt foundation.
He also thanks Edward Dahl from DWave for suggesting the problem of computing the dual
description of $\CORP(K_{4,4})$.

\bibliographystyle{amsplain_initials_eprint}
\bibliography{LatticeRef}

\providecommand{\bysame}{\leavevmode\hbox to3em{\hrulefill}\thinspace}
\providecommand{\MR}{\relax\ifhmode\unskip\space\fi MR }
\providecommand{\MRhref}[2]{%
  \href{http://www.ams.org/mathscinet-getitem?mr=#1}{#2}
}
\providecommand{\href}[2]{#2}
\begin{thebibliography}{10}

\bibitem{PolymakeLin}
B.~Assarf, E.~Gawrilow, K.~Herr, M.~Joswig, B.~Lorenz, A.~Paffenholz, and
  T.~Rehn, \emph{polymake in linear and integer programming}, preprint at {\tt
  arxiv:arXiv:1408.4653}, August 2014.

\bibitem{LeggettGargIneq}
D.~Avis, P.~Hayden, and M.~M. Wilde, \emph{Leggett-{G}arg inequalities and the
  geometry of the cut polytope}, Phys. Rev. A (3) \textbf{82} (2010), no.~3,
  030102, 4, URL: \url{http://dx.doi.org/10.1103/PhysRevA.82.030102},
  doi:10.1103/PhysRevA.82.030102.

\bibitem{AvisCorr}
D.~Avis, H.~Imai, and T.~Ito, \emph{On the relationship between convex bodies
  related to correlation experiments with dichotomic observables}, J. Phys. A
  \textbf{39} (2006), no.~36, 11283--11299, URL:
  \url{http://dx.doi.org/10.1088/0305-4470/39/36/010},
  doi:10.1088/0305-4470/39/36/010.

\bibitem{AvisTwoParty}
D.~Avis, H.~Imai, T.~Ito, and Y.~Sasaki, \emph{Two-party {B}ell inequalities
  derived from combinatorics via triangular elimination}, J. Phys. A
  \textbf{38} (2005), no.~50, 10971--10987, URL:
  \url{http://dx.doi.org/10.1088/0305-4470/38/50/007},
  doi:10.1088/0305-4470/38/50/007.

\bibitem{OrbitFacetCutPolytope6}
D.~Avis and Mutt, \emph{All the facets of the six-point {H}amming cone},
  European J. Combin. \textbf{10} (1989), no.~4, 309--312, URL:
  \url{http://dx.doi.org/10.1016/S0195-6698(89)80002-2},
  doi:10.1016/S0195-6698(89)80002-2.

\bibitem{BalinskiTh}
M.~L. Balinski, \emph{On the graph structure of convex polyhedra in
  {$n$}-space}, Pacific J. Math. \textbf{11} (1961), 431--434.

\bibitem{BarahonaCutContract}
F.~Barahona, \emph{The max-cut problem on graphs not contractible to
  {$K_{s}$}}, Oper. Res. Lett. \textbf{2} (1983), no.~3, 107--111, URL:
  \url{http://dx.doi.org/10.1016/0167-6377(83)90016-0},
  doi:10.1016/0167-6377(83)90016-0.

\bibitem{BarahonaMahjoubOnCutPolytope}
F.~Barahona and A.~R. Mahjoub, \emph{On the cut polytope}, Math. Programming
  \textbf{36} (1986), no.~2, 157--173, URL:
  \url{http://dx.doi.org/10.1007/BF02592023}, doi:10.1007/BF02592023.

\bibitem{LMS_JCM_groupPolytope}
D.~Bremner, M.~Dutour~Sikiri{\'c}, D.~V. Pasechnik, T.~Rehn, and
  A.~Sch{\"u}rmann, \emph{Computing symmetry groups of polyhedra}, LMS J.
  Comput. Math. \textbf{17} (2014), no.~1, 565--581.

\bibitem{symsurvey}
D.~Bremner, M.~Dutour~Sikiri{\'c}, and A.~Sch{\"u}rmann, \emph{Polyhedral
  representation conversion up to symmetries}, Polyhedral computation, CRM
  Proc. Lecture Notes, vol.~48, Amer. Math. Soc., Providence, RI, 2009,
  pp.~45--71.

\bibitem{CR}
T.~Christof and G.~Reinelt, \emph{Decomposition and parallelization techniques
  for enumerating the facets of combinatorial polytopes}, Internat. J. Comput.
  Geom. Appl. \textbf{11} (2001), no.~4, 423--437, URL:
  \url{http://dx.doi.org/10.1142/S0218195901000560},
  doi:10.1142/S0218195901000560.

\bibitem{DualCutPolytope}
A.~Deza and M.~Deza, \emph{On the skeleton of the dual cut polytope}, Jerusalem
  combinatorics '93, Contemp. Math., vol. 178, Amer. Math. Soc., Providence,
  RI, 1994, pp.~101--111, URL: \url{http://dx.doi.org/10.1090/conm/178/01910},
  doi:10.1090/conm/178/01910.

\bibitem{OrbitVerticesMetricPolytope7}
A.~Deza, M.~Deza, and K.~Fukuda, \emph{On skeletons, diameters and volumes of
  metric polyhedra}, Combinatorics and computer science ({B}rest, 1995),
  Lecture Notes in Comput. Sci., vol. 1120, Springer, Berlin, 1996,
  pp.~112--128, URL: \url{http://dx.doi.org/10.1007/3-540-61576-8_78},
  doi:10.1007/3-540-61576-8-78.

\bibitem{OrbitPolytopeMetricPolytope8}
A.~Deza, K.~Fukuda, T.~Mizutani, and C.~Vo, \emph{On the face lattice of the
  metric polytope}, Discrete and computational geometry, Lecture Notes in
  Comput. Sci., vol. 2866, Springer, Berlin, 2003, pp.~118--128, URL:
  \url{http://dx.doi.org/10.1007/978-3-540-44400-8_12},
  doi:10.1007/978-3-540-44400-8-12.

\bibitem{Relatives}
M.~Deza, V.~P. Grishukhin, and M.~Laurent, \emph{The symmetries of the cut
  polytope and of some relatives}, Applied geometry and discrete mathematics,
  DIMACS Ser. Discrete Math. Theoret. Comput. Sci., vol.~4, Amer. Math. Soc.,
  Providence, RI, 1991, pp.~205--220.

\bibitem{Hyp8}
M.~Deza and M.~Dutour~Sikiri\'c, \emph{{The facets and extreme rays of the
  hypermetric cone on eight vertices}}, in preparation.

\bibitem{DL}
M.~M. Deza and M.~Laurent, \emph{Geometry of cuts and metrics}, Algorithms and
  Combinatorics, vol.~15, Springer, Heidelberg, 2010, First softcover printing
  of the 1997 original [MR1460488], URL:
  \url{http://dx.doi.org/10.1007/978-3-642-04295-9},
  doi:10.1007/978-3-642-04295-9.

\bibitem{WebPageCutPolytopes}
M.~Dutour~Sikiri\'c, \emph{Cut polytopes}, URL:
  \url{http://mathieudutour.altervista.org/CutPolytopes/}.

\bibitem{Polyhedral}
M.~Dutour~Sikiri\'c, \emph{Polyhedral}, URL:
  \url{http://mathieudutour.altervista.org/Polyhedral/}.

\bibitem{PerfectDim8}
M.~Dutour~Sikiri{\'c}, A.~Sch{\"u}rmann, and F.~Vallentin, \emph{Classification
  of eight-dimensional perfect forms}, Electron. Res. Announc. Amer. Math. Soc.
  \textbf{13} (2007), 21--32 (electronic), URL:
  \url{http://dx.doi.org/10.1090/S1079-6762-07-00171-0},
  doi:10.1090/S1079-6762-07-00171-0.

\bibitem{ComplexityVoronoiDSV}
M.~Dutour~Sikiri{\'c}, A.~Sch{\"u}rmann, and F.~Vallentin, \emph{Complexity and
  algorithms for computing {V}oronoi cells of lattices}, Math. Comp.
  \textbf{78} (2009), no.~267, 1713--1731, URL:
  \url{http://dx.doi.org/10.1090/S0025-5718-09-02224-8},
  doi:10.1090/S0025-5718-09-02224-8.

\bibitem{ContactLeech}
M.~Dutour~Sikiri{\'c}, A.~Sch{\"u}rmann, and F.~Vallentin, \emph{The contact
  polytope of the {L}eech lattice}, Discrete Comput. Geom. \textbf{44} (2010),
  no.~4, 904--911, URL: \url{http://dx.doi.org/10.1007/s00454-010-9266-z},
  doi:10.1007/s00454-010-9266-z.

\bibitem{FineHidden}
A.~Fine, \emph{Hidden variables, joint probability, and the {B}ell
  inequalities}, Phys. Rev. Lett. \textbf{48} (1982), no.~5, 291--295, URL:
  \url{http://dx.doi.org/10.1103/PhysRevLett.48.291},
  doi:10.1103/PhysRevLett.48.291.

\bibitem{OrbitFacetCutPolytope7}
V.~P. Grishukhin, \emph{All facets of the cut cone {${\bf C}_n$} for {$n=7$}
  are known}, European J. Combin. \textbf{11} (1990), no.~2, 115--117, URL:
  \url{http://dx.doi.org/10.1016/S0195-6698(13)80064-9},
  doi:10.1016/S0195-6698(13)80064-9.

\bibitem{OrbitRaysMetricCone7}
V.~P. Grishukhin, \emph{Computing extreme rays of the metric cone for seven
  points}, European J. Combin. \textbf{13} (1992), no.~3, 153--165, URL:
  \url{http://dx.doi.org/10.1016/0195-6698(92)90021-Q},
  doi:10.1016/0195-6698(92)90021-Q.

\bibitem{PitowskyOptimal}
I.~Pitowsky and K.~Svozil, \emph{Optimal tests of quantum nonlocality}, Phys.
  Rev. A (3) \textbf{64} (2001), no.~1, 014102, 4, URL:
  \url{http://dx.doi.org/10.1103/PhysRevA.64.014102},
  doi:10.1103/PhysRevA.64.014102.

\bibitem{SeymourMatroidMulticommodity}
P.~D. Seymour, \emph{Matroids and multicommodity flows}, European J. Combin.
  \textbf{2} (1981), no.~3, 257--290, URL:
  \url{http://dx.doi.org/10.1016/S0195-6698(81)80033-9},
  doi:10.1016/S0195-6698(81)80033-9.

\bibitem{OrbitFacetCutPolytope5}
M.~E. Tylkin (=M.~Deza), \emph{On {H}amming geometry of unitary cubes}, Soviet
  Physics. Dokl. \textbf{5} (1960), 940--943.

\bibitem{WernerBellInequality}
R.~Werner, \emph{All the bell inequalities}, URL:
  \url{http://qig.itp.uni-hannover.de/qiproblems/All_the_Bell_Inequalities}.

\end{thebibliography}

\end{document}